\documentclass[leqno,12pt]{article}
\usepackage{amsmath,amssymb,latexsym}
\usepackage{a4wide,theorem}
\def\Rad{\mathop{\mathrm{Rad}}\nolimits}

\def\min{\mathop{\mathrm{min}}\nolimits}

\theoremstyle{break}
\newtheorem{sahilf}{Proposition}[section]
\newtheorem{kohilf}[sahilf]{Corollary}
\newtheorem{thhilf}[sahilf]{Theorem}
\newtheorem{lehilf}[sahilf]{Lemma}
\newtheorem{fact}[sahilf]{Fact}
\theoremstyle{plain}
\theorembodyfont{\rmfamily}
\newtheorem{dehilf}[sahilf]{Definition}
\newtheorem{exhilf}[sahilf]{Example}
\newtheorem{bem}[sahilf]{Remark}

\def\qedbox{\mbox{\large $\Box$}}

\newcommand{\deflabel}[1]{\bf1\hfill}%
{\begin{list}{}%
{\settowidth{\labelwidth}{\bf1}%
\setlength{\leftmargin}{\labelwidth}%
\addtolength{\leftmargin}{\labelsep}%
}}%
{\end{list}}%

\newcommand{\gen}[1]{\langle #1 \rangle}

\newenvironment{Bw}{\begin{list}{}%
       {\setlength{\leftmargin}{0mm}}%
       \item {\bf Proof:}
             }{\rule{0mm}{0mm}\nolinebreak\hfill\qedbox
                \end{list}}

\newenvironment{sketch}{\begin{list}{}%
       {\setlength{\leftmargin}{0mm}}%
       \item {\bf Sketch of proof:}
             }{\rule{0mm}{0mm}\nolinebreak\hfill\qedbox
                \end{list}}

\begin{document}

\title{The sphericity of the complex of non-degenerate subspaces} 

\author{Alice Devillers\thanks{The first author was a Collaborateur Scientifique F.R.S.-FNRS during the time when most of the work for this article was done.}$\ \ ^1$
, Ralf Gramlich\thanks{The second author gratefully acknowledges a Heisenberg fellowship of the Deutsche Forschungsgemeinschaft.}$\ \ {}^{2,3}$, and Bernhard M\"uhlherr$^4$} 

\maketitle

{\footnotesize
\noindent 1: University of Western Australia,
School of Mathematics and Statistics,
35 Stirling Highway, Crawley 6009, Western Australia, e-mail: \texttt{adevil@maths.uwa.edu.au} \\
2: TU Darmstadt, FB Mathematik, Schlossgartenstrasse 7, 64289 Darmstadt, Germany, \\ e-mail: \texttt{gramlich@mathematik.tu-darmstadt.de} \\
3: The University of Birmingham, School of Mathematics, Edgbaston, Birmingham, B15 2TT, United Kingdom, e-mail: \texttt{ralfg@maths.bham.ac.uk} \\
4: Universit\"at Gie\ss en, Mathematisches Institut, Arndtstrasse 2, 35392 Gie\ss en, Germany  \\ e-mail: \texttt{muhlherr@math.uni-giessen.de}
}

\begin{abstract} 
\noindent We prove that the complex of proper non-trivial non-degenerate subspaces of a finite-dimensional vector space endowed with a non-degenerate sesquilinear form is homotopy equivalent to a wedge of spheres. Additionally, we show that the same is true for a slightly wider class of simplicial complexes, the so-called generalized Phan geometries of type $A_n$. These generalized Phan geometries occur as relative links of the filtration studied in \cite{DM}, whose sphericity implies topological finiteness properties of suitable arithmetic groups and allows for a revision of Phan's group-theoretical local recognition \cite{Phan:1977}, \cite{Phan:1977a} of suitable finite groups of Lie type with simply laced diagrams.  
\end{abstract}

\section{Introduction}

In this paper we prove the following result.

\medskip
\noindent {\bf Main Theorem} \\
{\em Let $\mathbb{K}$ be an arbitrary field, 
let $\sigma$ be an automorphism of $\mathbb{K}$ of order one or two and 
let $V$ be an $(n+1)$-dimensional vector space over $\mathbb{K}$ endowed 
with a non-degenerate $\sigma$-hermitian form $h$. Moreover, 
let $\Gamma$ be the simplicial complex whose vertices are 
the non-trivial subspaces which are non-degenerate with
respect to $h$ and where incidence is defined by inclusion.  
If $\mathbb{K}$ is a finite field of order $q$, assume $2^n < q$ in case $\sigma = \mathrm{id}$ and $2^{n-1}(\sqrt{q}+1) < q$ in case $\sigma \neq \mathrm{id}$. 

Then $\Gamma$ is Cohen-Macaulay. In particular, $|\Gamma|$ is homotopy equivalent to a wedge of $(n-1)$-spheres.
}

\medskip
For $\mathbb{K} = \mathbb{F}_{q^2}$, $q \geq 4$, $n \geq 3$, and $h$ hermitian, the above simplicial complex $\Gamma$ is known to be simply connected by \cite{Aschbacher}, \cite{Das:1994}, \cite{Bennett/Shpectorov:2004}. Using standard arguments from geometric group theory this implies Phan's theorem \cite{Phan:1977}, which states that the group $\mathrm{SU}_{n+1}(q^2)$ equals the universal enveloping group of the amalgam of its fundamental subgroups of rank one and two. Because of this prominent role in the revision of Phan's results, the geometry of non-degenerate subspaces described above is called
{\it the Phan geometry associated to the pair $(V,h)$}.

In fact, there is
a canonical way to associate a Phan geometry to any flip 
of a spherical building  
(see \cite{Bennett/Gramlich/Hoffman/Shpectorov:2003}). 
A flip is an involutory automorphism of the spherical building that acts via conjugation with the longest word of the Weyl group on the Weyl distance and maps some chamber to an opposite chamber. In other words, a flip is an involutory isometry of the associated spherical twin building (see \cite{Tits:1992}) such that there exists a chamber that is opposite its image under that involution. 
It was observed in \cite{DM} 
that a Phan geometry induced by a flip shares a lot of properties 
with the geometry far away from a chamber. In a sense, our
work underlines this observation, because we successfully adapted the strategy of proof of the analogue of the Main Theorem for the far-away geometries in $A_n$ buildings
obtained by Abels and Abramenko in \cite{Abels/Abramenko:1993} (see also \cite{Abramenko:1996}) to the situation of Phan geometries.

The sphericity of the far-away geometries turns out to be
important in the work by Abels and Abramenko on topological finiteness properties of the group $\mathrm{SL}_{n+1}(\mathbb{F}_q[t])$.
The main motivation for our investigations was to obtain
a similar tool for proving topological finiteness properties
for centralizers of flips in
split Kac-Moody groups over finite fields, such as for 
the group      
 $$\mathrm{SU}_{n+1}(\mathbb{F}_{q^2}[t,t^{-1}], \theta) = 
\{ A \in \mathrm{SL}_{n+1}(\mathbb{F}_{q^2}[t,t^{-1}]) 
\mid (A^{-1})^\theta = A^T \}$$
where $\theta$ is the involution 
of the field $\mathbb{F}_{q^2}(t)$ mapping $t$ to $t^{-1}$ and acting 
as the Frobenius involution on $\mathbb{F}_{q^2}$.
However, it turned out that the statement of the Main Theorem
is not sufficient for this purpose. The reason 
is that the relative links in the filtration of \cite{DM}
are slightly more general than the Phan geometries.
This lead us to the definition of {\it generalized Phan geometries}, see
Section \ref{def}, which is quite technical. Fortunately, while modifying 
an earlier version of the proof of the Main Theorem in order
to  make it work for 
the generalized Phan geometries, we observed that the arguments
became much more natural in this wider context and some technicalities just
disappeared. In fact, the class of generalized Phan geometries contains both the class of Phan geometries and the class of geometries far from a chamber. 

As we have already indicated, the result
of the present paper is meant to provide a tool to determine 
the finiteness lengths
of unitary forms of split Kac-Moody groups over finite fields
along the lines of \cite{Abramenko:1996}. Apart from this application
we would like to mention that the cases of small dimension $n \in \{ 2, 3 \}$ of our result are
precisely what is needed in the revision of Phan's result 
for the finite Chevalley groups with simply laced diagram (see \cite{GHMS}; see \cite{Phan:1977a} for Phan's original result).
Our proof provides a more systematic approach
than the one given in \cite{GHMS}, which is based on a
case by case distinction. However, the proof in \cite{GHMS} covers $q=4$, which the proof presented in this paper does not.

The methods needed in order to obtain the Main Theorem are quite elementary and
consist of linear algebra and basic homology
 theory. Our approach is
based on a combination of the
strategy of the sphericity proof for far-away geometries
given in \cite{Abramenko:1996} and several techniques used in the
revision of Phan's
presentation \cite{Bennett/Gramlich/Hoffman/Shpectorov:2003}.
The necessary material from topology is 
discussed in Section \ref{basics}, while in Section \ref{def} 
we study some examples of generalized Phan geometries of type $A_n$. 
Section \ref{mainresult} concentrates on proving that  
generalized Phan geometries of type $A_n$ are spherical over sufficiently large fields. 
This goal is achieved in Theorem \ref{allrad} and its Corollary \ref{allrad2}. The Main Theorem stated above
is just a special case.
Finally, in Section \ref{outline} we outline in more detail 
the two applications of this 
sphericity result to group theory mentioned above. 
While Sections \ref{basics}, \ref{def} and \ref{mainresult} are 
largely self-contained, 
the character of Section \ref{outline} is that of a survey. 
The final section is intended to serve as a motivation for our interest 
in generalized Phan 
geometries of type $A_n$, 
so that our arguments given there should only be considered as sketches.
The content of Section \ref{outline} is based on Fact \ref{fact}, which has already been proved in special situations in the literature, cf.\ \cite{Blok/Hoffman}. We intend to a give a general building-theoretic proof in a separate publication.

\medskip
\noindent {\bf Acknowledgement:} The research resulting in this paper has been 
initiated during an RiP stay at Oberwolfach in summer 2006. 
The authors express their gratitude for the hospitality of the 
Mathematisches Forschungsinstitut Oberwolfach. The authors thank Peter 
Abramenko, Kai-Uwe Bux, Max Horn, William Kantor, Linus Kramer, Sergey 
Shpectorov, and Stefan Witzel for several helpful remarks and comments. The
authors also thank the anonymous referee for an excellent report that lead 
to a considerable improvement of this
article.

\section{Basics} \label{basics}

\begin{dehilf}[Simplicial complex]
A {\em} simplicial complex $K$ consists of a set $V$ of vertices and a set $\mathfrak{S}$ of finite non-empty subsets of $V$ called {\em simplices} such that
\begin{enumerate}
\item any set consisting of exactly one vertex is a simplex, and
\item any non-empty subset of a simplex is a simplex.
\end{enumerate}
A simplicial complex $K$ is called {\em pure} if each simplex is contained in a maximal simplex and if, moreover, all maximal simplices have the same cardinality $n$; in this case we say that the simplicial complex $K$ has {\em dimension} $n-1$.
The {\em join} of two simplicial complexes $K_1=(V_1, \mathfrak{S}_1)$ and $K_2=(V_2, \mathfrak{S}_2)$, in symbols $K_1 * K_2$, equals the simplicial complex with vertex set equal to the disjoint union $V_1 \sqcup V_2$ and set of simplices equal to $\{ v \cup w \in V_1 \sqcup V_2 \mid  v \in \mathfrak{S}_1, w \in \mathfrak{S}_2 \}$.     

Furthermore, let $|K|$, called {\em realization of $K$}, be the set of all functions $\alpha$ from the set of vertices of $K$ to the real unit interval $I$ such that, for any $\alpha$, the set $\{ v \in K \mid \alpha(v) \neq 0 \}$ is a simplex of $K$, and $\sum_{v \in K} \alpha(v) = 1$, i.e., $|K|$ is obtained from $K$ via barycentric coordinates.
\end{dehilf}

In this article we consider the weak (coherent) topology on $|K|$, cf.\ \cite[3.1.14]{Spanier}. By \cite[7.6.4]{Spanier} the realization $|K|$ of a simplicial complex $K$ is a CW complex.

\begin{dehilf}[Realization]
A poset $X$ gives rise to a simplicial complex  with the non-empty finite chains of $X$ as simplices. Endowed with the weak topology, this simplicial complex is called the {\em realization} of $X$ and denoted by $|X|$. 

In analogy to simplicial complexes, the {\em join} of two posets $X$ and $Y$, also denoted by $X * Y$, is the disjoint union $X \sqcup Y$ endowed with the partial order which agrees with the given partial orders on $X$ and on $Y$ and which is such that any element of $X$ is less than any element of $Y$.
\end{dehilf}

\begin{dehilf}[Closed simplex, open simplex, star, link]
Let $K$ be a simplicial complex and let $s$ be a simplex of $K$. The {\em closed simplex} $|s|$ equals $\{ \alpha \in |K| \mid \alpha(v) \neq 0 \Longrightarrow v \in s \}$. Hence, if $s$ is a $q$-simplex, then $|s|$ is naturally homeomorphic to the standard simplex $\{ x \in \mathbb{R}^{q+1} \mid 0 \leq x_i \leq 1, \sum x_i = 1 \}$. The {\em open simplex} $\gen{s}$ equals $\{ \alpha \in |K| \mid \alpha(v) \neq 0 \Longleftrightarrow v \in s \}$. The {\em open simplex} $\gen{s}$ is open in $|s|$, but not necessarily in $|K|$.

For a vertex $v \in K$, the (open) {\em star} $St_K(v)$ equals $\{ \alpha \in |K| \mid \alpha(v) \neq 0 \}$.
The {\em link} of $x$ is defined as $Lk_K(x) := \overline{St_K(x)} \backslash St_K(x)$.
\end{dehilf}

We have $St_K(x) = \bigcup \{ \gen{s} \mid \mbox{$x$ is vertex of $s$} \}$, cf.\ \cite[3.1.24]{Spanier}.

\begin{dehilf}[Wedge sum]
Let $X$ and $Y$ be pointed spaces, i.e., topological spaces with distinguished base points $x_0$ and $y_0$. Then the {\em wedge sum} $X \vee Y$ of $X$ and $Y$ is the quotient of the disjoint union $X \sqcup Y$ by the identification $x_0 \sim y_0$, i.e., $$X \vee Y := (X \sqcup Y) / \{ x_0 \sim y_0 \}.$$
In general, if $(X_i)_{i \in I}$ is a family of pointed spaces with base points $(x_i)_{i \in I}$ and if $x$ is a one-point space, then the {\em wedge sum} of this family is given by $$\bigvee_{i \in I} X_i := \bigsqcup_{i \in I} X_i / \{ x_i \sim x \}.$$ 
The wedge sum of a family of spheres each of which is homeomorphic to the standard sphere $\mathbb{S}^n = \{ x \in \mathbb{R}^{n+1} \mid \| x \| = 1 \}$ is called a {\em wedge of spheres} or, if one wants to specify the dimension, a {\em wedge of $n$-spheres}.
\end{dehilf}

The following topological observations are crucial for the proof of our main result. 
A CW complex of dimension at most $m$ is called $m$-\emph{spherical} if it 
has the homotopy
type of a wedge $m$-spheres.

\begin{lehilf} \label{above}
A connected $1$-dimensional CW complex $X$ has the homotopy type of a wedge of
$1$-spheres.
\end{lehilf}

\begin{Bw}
Let $T\subseteq X$ be a maximal subtree. Then $T$ is contractible
and contains all vertices (see \cite[Lemmas 3.7.1 and 3.7.2]{Spanier}), whence by the homotopy extension property $X$ is homotopy equivalent to $X/T$ (see \cite[Corollary 3.2.5]{Spanier}). This quotient is a connected
$1$-dimensional
complex with exactly one vertex, i.e. a wedge of $1$-spheres.
\end{Bw}

\begin{lehilf} \label{mainlemma}
Let $X$ be a CW complex of dimension at most $m$ which is a union
$X=B\cup\bigcup_iA_i$ of subcomplexes $B,A_i$, where $i$ ranges over some
(possibly infinite) index set $I$. Suppose that $B$ is $m$-spherical, that 
the
$A_i$ are contractible, that $A_i\cap A_j\subseteq B$ for $i\neq j$, and 
that
$A_i\cap B$ is $m-1$-spherical, for every $i$. Then $X$ is $m$-spherical.
\end{lehilf}

\begin{Bw}
If $X$ is $1$-dimensional, the 
result is true by Lemma \ref{above}.

Assume now that $m\geq 2$. 
Let $Y=B\cup A_i$. The Mayer-Vietoris sequence in reduced homology is
\[
 0 \to \tilde H_m(B)\oplus \tilde H_m(A_i)\to \tilde H_m(Y)\to 
\tilde H_{m-1}(B\cap A_i)\to 0
\]
so $Y$ has the same homology as a wedge of $m$-spheres. By Van Kampen's
Theorem, $B\cap A_i$ is $1$-connected. It follows from Hurewicz' Theorem 
\cite[7.5.4]{Spanier} that
$\pi_*(Y)\cong H_*(Y)$ for $* \leq m$. Let $f:\bigvee_JS^m\to Y$ be a map sending each
sphere summand to a free generator of $\pi_m(Y)$. Then $f$ is a homology
equivalence and therefore a homotopy equivalence \cite[7.5.9, 7.6.24]{Spanier}.
It follows by induction that $B\cup A_{i_1}\cup\cdots\cup A_{i_r}$
has the homotopy type of a wedge of $m$-spheres.
Since every compact set $C\subseteq X$ is contained in such a finite
union, $X$ is $m-1$-connected and has the same homology as a wedge of 
$m$-spheres.
By the same argument as before, we find a homotopy equivalence
$\bigvee_KS^m\to X$.
\end{Bw}

We close this section by specifying the characterizing property of Cohen-Macaulayness that we are going to use in this article; cf.\ \cite[Section 8]{Quillen:1978}.

\begin{dehilf}[Cohen-Macaulay]
A $d$-dimensional simplicial complex $K$ is called {\em Cohen-Macaulay}, if it is $d$-spherical, and if in addition the link of each $p$-simplex is $(d-p-1)$-spherical. 
%
%
%
%
%for each simplex $s$ of $K$ the reduced homology group $\tilde H_p(Lk_K(s))$ is trivial for all $p < \mathrm{dim}(Lk_K(s))$.
\end{dehilf}

\section{Generalized Phan geometries of type $A_n$} \label{def}

\begin{dehilf} \label{trans}
 Let $V$ be a vector space. Two subspaces $A$ and $B$ of $V$ are {\em opposite} when $V=A\oplus B$. A flag is a chain of incident subspaces. A subspace $A$ is {\em transversal} or {\em in general position} to a flag $F$, denoted by $A\pitchfork_V F$, if for any subspace $B$ of $F$ we have $A\cap B=\{0\}$ or $V=A+B$. 
\end{dehilf}
Notice that $A\pitchfork_V F$ if and only if there is a subspace $C$ of $V$ incident with $F$ such that $A$ and $C$ are opposite.

\begin{dehilf}
Let $\mathbb{K}$ be a field, let $\sigma$ be an automorphism of $\mathbb{K}$  of order one or two, let $V$ be an $(n+1)$-dimensional $\mathbb{K}$-vector space, let $0 \leq t \leq n$, let $F = \{ V_i \mid 0 \leq i \leq t+1\}$ with $0 = V_0 \lneq V_1 \lneq \cdots \lneq V_t \lneq V_{t+1} = V$ be a flag of subspaces of $V$, and let $\omega_{i} : V_{i+1} \times V_{i+1} \to \mathbb{K}$, $0 \leq i \leq t$, be $\sigma$-hermitian forms which admit non-isotropic vectors and satisfy $$\Rad\left(\omega_i\right) = V_i.$$ 
The {\em generalized Phan geometry of type $A_n$} of $V$ with respect to $F$ and $(\omega_{i})_{0 \leq i \leq t}$ consists of all proper non-trivial vector subspaces $U$ of $V$ transversal to $F$ with $U \cap V_{k_U+1}$ non-degenerate with respect to $\omega_{k_U}$ where $k_U = \min \{ i \in \{0, \ldots, t \} \mid U \cap V_{i+1} \not= \{0\} \}$. 
\end{dehilf}

\begin{exhilf}\label{examples}
For $t = 0$ the generalized Phan geometry on $V$ equals the complex of proper non-trivial non-degenerate subspaces of $V$ with respect to a non-degenerate $\sigma$-hermitian form. For $\sigma$ of order one and the field $\mathbb{K}$ of characteristic distinct from two, this complex is simply connected by \cite{Altmann:2003}, \cite{Altmann/Gramlich:2006}. For $\mathbb{K} = \mathbb{F}_{q^2}$ and $\sigma$ the Frobenius involution, this complex has been studied in \cite{Aschbacher}, \cite{Das:1994} and in \cite{Bennett/Shpectorov:2004} in order to re-prove Phan's characterization \cite{Phan:1977} of the group $\mathrm{SU}_{n+1}(q^2)$ via combinatorial topology. 
\end{exhilf}

\begin{exhilf} \label{examples2}
 For $t = n$, the generalized Phan geometry on $V$ equals the geometry opposite a chamber in a building of type $A_n$ (whose sphericity has been established in \cite{Abels/Abramenko:1993}, see also \cite{Abramenko:1996}).
Indeed, the flag $F$ is a chamber and each ${\omega_i}$ has rank one with radical $V_i$. Therefore any vector $v \in V_{i+1} \backslash V_i$ is non-degenerate with respect to $\omega_i$, so that any subspace $U$ of $V$ with $U \oplus V_i = V$ intersects $V_{i+1}$ in a non-degenerate (with respect to $\omega_i$) one-dimensional subspace. 
\end{exhilf}

Just like the geometries opposite a chamber occur as relative links in filtrations used to establish finiteness properties of certain arithmetic groups (cf.\ \cite{Abels:1991}, \cite{Abramenko:1996}), so do the generalized Phan geometries. The following series of examples illustrates how generalized Phan geometries occur as relative links.

\begin{exhilf}
Let $V$ be a vector space and consider a generalized Phan geometry on $V$ with $t=0$, i.e., the geometry of non-degenerate subspaces of $V$ with respect to a form $\omega_0 : V \times V \to \mathbb{K}$. Let $H$ be a hyperplane of $V$ with radical $R$ of dimension zero or one. Then the intersection of the link of $H$ (in the complex of proper non-trivial subspaces of $V$) with the subcomplex of proper non-trivial $\omega_0$-non-degenerate subspaces of $V$, equals all proper non-trivial subspaces of $H$ which are non-degenerate with respect to ${\omega_0}_{|H \times H}$. Each of these subspaces is transversal to the flag $\{ \{ 0 \}, R, H \}$, so that this relative link is a generalized Phan geometry. (Notice that the form on $R$ does not play any role for this generalized Phan geometry.) 
\end{exhilf}

In Lemma \ref{residue} we prove that the relative link of an arbitrary element of a generalized Phan geometry is isomorphic to a generalized Phan geometry or the direct sum (i.e., join) of two generalized Phan geometries.

\begin{exhilf} \label{examples3}
Let $\mathcal{B}$ be a building of type $A_n$, let $R$ be a residue of $\mathcal{B}$, and let $\pi$ be a polarity of $\mathcal{B}$. Using the notation from \cite{DM}, the chamber system $A_\pi(R)$ consists of all those chambers of $R$ with maximal possible distance (i.e.,  minimal possible codistance) from their image under $\pi$. This chamber system plays a crucial role in the filtration studied in \cite{DM} and yields a direct sum of generalized Phan geometries, cf.\ Fact \ref{fact}. When dealing with polarities of projective spaces (buildings of type $A_n$), this distance of a chamber $c$ from its image under $\pi$ can be measured by the grade of nondegeneracy of the elements contained in the chamber $c$. 

For instance, let $V$ be a six-dimensional vector space (which corresponds to a building of type $A_5$) endowed with a non-degenerate unitary form $h$ (which corresponds to a polarity $\pi$). Consider the residue of $V$ with respect to the flag consisting of an $h$-singular one-dimensional subspace $p$ and an incident five-dimensional subspace $H$ with one-dimensional $h$-radical $x$. Fixing a four-dimensional complement $C$ of $p$ not containing $x$, we can identify the set of chambers of $V$ containing $p$ and $H$ with the set of chambers of $C$ by taking the intersection. The set of chambers in the residue of $(p,H)$ which are mapped by $\pi$ to a chamber as far away as possible, correspond in $C$ to maximal flags of subspaces of the following types: (i) one-dimensional subspaces of $C$ opposite to $p^\perp \cap C$, (ii) two-dimensional subspaces of $C$ transversal to $p^\perp \cap C$ which intersect $p^\perp \cap C$ in an $h$-non-degenerate one-dimensional subspace, and (iii) three-dimensional subspaces of $C$ transversal to $p^\perp \cap C$ which intersect $p^\perp \cap C$ in an $h$-non-degenerate two-dimensional subspace. This description shows that this residue is isomorphic to a generalized Phan geometry of $C$ with respect to the flag $\{ \{ 0 \}, C \cap p^\perp, C \}$ and the form $h_{|(C \cap p^\perp) \times (C \cap p^\perp)}$ and an arbitrary sesquilinear form $f : C \times C \to \mathbb{K}$ of rank one with radical $C \cap p^\perp$.
\end{exhilf}

%\begin{exhilf}
%For $n=1$, a generalized Phan geometry is the set of non-degenerate one-dimensional subspaces of a two-dimensional vector space endowed with a non-trivial hermitian form, the form being non-degenerate ($t=0$) or with radical of dimension one ($t=1$).
%\end{exhilf}

%\begin{exhilf}
%For $n=2, 3$ and $\mathbb{K} = \mathbb{F}_{q^2}$, the generalized Phan geometries are exactly the residues considered in \cite{DM} in order to prove simple connectedness of the Phan geometries of type $\Delta$ over $\mathbb{F}_{q^2}$ where $\Delta$ is a simply laced three-spherical diagram.
%\end{exhilf}

\section{Sphericity of generalized Phan geometries} \label{mainresult}

In this section we prove the Main Theorem from the introduction in the guise of Theorem \ref{allrad} and Corollary \ref{allrad2}. We achieve this by using a filtration method similar to the one described in \cite[Theorem 1.1]{Abels/Abramenko:1993}, \cite[Chapter II, Proposition 12]{Abramenko:1996}. Our incarnation of that strategy, as described in the proof of Theorem \ref{allrad} below, requires that there is a canonical way to extend each form $\omega_i$ of a generalized Phan geometry to the whole ambient vector space. The following lemma takes care of this issue.

\begin{lehilf} \label{extension}
Let $\Gamma$ be a generalized Phan geometry of $V$ with respect to a flag $F = \{ V_i \mid 0 \leq i \leq t+1 \}$ with forms $(\omega_i)_{0 \leq i \leq t}$. Let $p$ be a one-dimensional space in $V$ that is $\omega_t$-non-degenerate and satisfies $p \cap V_t = \{ 0 \}$. Then, for $0 \leq i \leq t$, there exist $\sigma$-hermitian forms $\underline{\omega}_i : V \times V \to \mathbb{K}$ with the properties
\begin{itemize}
\item ${\underline{\omega}_i}_{|V_{i+1} \times V_{i+1}} = \omega_i$,
\item $\Rad\left(\underline{\omega}_i\right) = \Rad\left({\omega_i}\right) = V_i$,
\item $p$ is $\underline{\omega}_i$-non-degenerate ($0 \leq i \leq t$) and satisfies $p^{\perp_{\underline{\omega}_i}} = p^{\perp_{\underline{\omega}_j}}$ for all $0 \leq i, j \leq t$.
\end{itemize}
\end{lehilf}

\begin{Bw}
Let $p$ be a $\omega_t$-non-degenerate one-dimensional subspace of $V$, which necessarily satisfies $p \cap V_t = \{ 0 \}$. Then 
\begin{eqnarray*}
V & = & C_1 \oplus C_2 \oplus C_3 \oplus \cdots \oplus C_t \oplus C_{t+1} \oplus p \\
& \cong & V_1 \oplus (V_2/V_1) \oplus (V_3/V_2) \oplus \cdots \oplus (V_t/V_{t-1}) \oplus (p^{\perp_{\omega_t}} / V_t) \oplus p,
\end{eqnarray*}
where $C_i$ is a complement of $V_{i-1}$ in $V_{i}$ for $1 \leq i \leq t$ and $C_{t+1}$ is a complement of $V_t$ in  $p^{\perp_{\omega_t}}$.
The definition $\underline{\omega}_i : V \times V \to \mathbb{K}$ via $$\underline{\omega}_i(x,y) := \left(\sum_{i \leq j \leq t} \omega_j\left(\mathrm{pr}_{C_{j+1}}(x),\mathrm{pr}_{C_{j+1}}(y)\right)\right) + \omega_{t}\left(\mathrm{pr}_{p}(x),\mathrm{pr}_{p}(y)\right)$$ finishes the proof, where $\mathrm{pr}_A$ denotes the projection on $A$ in the above direct decomposition of $V$.
\end{Bw}

The purpose of the following lemma is two-fold. On one hand it allows us to make use of an induction on the dimension when studying the homotopy type of a generalized Phan geometry in Theorem \ref{allrad}. On the other hand it provides us with the local information required for the Cohen-Macaulayness of a generalized Phan geometry addressed in Corollary \ref{allrad2}.

\begin{lehilf} \label{residue}
Let $\Gamma$ be a generalized Phan geometry of $V$ with respect to a flag $F = \{ V_i \mid 0 \leq i \leq t+1 \}$ with forms $(\omega_{i})_{0 \leq i \leq t}$, and let $U \in \Gamma$ be an element of $\Gamma$. 

Then the geometries $\Gamma^{<U} := \{ x \in \Gamma \mid x \subsetneq U \}$ and $\Gamma^{>U} := \{ x \in \Gamma \mid x \supsetneq U \}$ are isomorphic to generalized Phan geometries of $U$, respectively, $V/U$. In particular, the link of $U$ in $|\Gamma|$ is isomorphic to the join $|\Gamma^{<U} * \Gamma^{>U}|$ of a generalized Phan geometry of $U$ and of a generalized Phan geometry of $V/U$.
\end{lehilf}

\begin{Bw}
%The residue $\Gamma_U$ clearly is the direct sum of $\Gamma^{<U} := \{ x \in \Gamma_U \mid x \leq U \}$ and $\Gamma^{>U} := \{ x \in \Gamma_U \mid x \geq U \}$, at most one of which is empty.
The first claim is immediate: The geometry $\Gamma^{<U}$ is isomorphic to the generalized Phan geometry of $U$ with respect to the flag $\{ V_i \cap U \mid k_U \leq i \leq t+1 \}$ and forms $\left({\omega_{i}}_{|(U \cap V_{i+1}) \times (U \cap V_{i+1})}\right)_{k_U \leq i \leq t}$ where $k_U = \min \{ i \in \{0, \ldots, t \} \mid U \cap V_{i+1} \not= \{0\} \}$. 

In order to establish the isomorphism type of $\Gamma^{>U}$, consider the vector space $\overline{V}=V/U$. Define $\overline{F}:=\{V_i + U \mid 0 \leq i \leq k_U+1\}$. This is a flag of $\overline{V}$ with $V_0 + U=0+U$ and $V_{k_U+1}+ U=\overline{V}$, so that any element of $\overline{V}$ can be written as $x+U$ with $x\in V_{k_U+1}$.
In fact, every element of $\overline{V}$ can be uniquely written as $x+U$ with $x \in (U \cap V_{k_U+1})^{\perp_{\omega_{k_U}}} =: W$. Notice that $V_{k_U}$ is contained in $W$, and therefore so are all subspaces $V_i$ for $0 \leq i \leq k_U$.
In this specified interval for $i$ define $\sigma$-hermitian forms  $$\overline{\omega}_i : (V_{i+1} + U) \times (V_{i+1} + U) \to \mathbb{K} : (x+U,y+U) \mapsto \overline{\omega}_i(x+U,y+U):=\omega_i(x,y)$$ 
where $x,y\in V_{i+1} \cap W$. With this definition, certainly $\Rad\left({\overline{\omega}_i}\right)= V_i + U$ for all $0 \leq i \leq k_U$, so that we have defined a generalized Phan geometry in $\overline{V}$.

It remains to prove that this generalized Phan geometry is isomorphic to $\Gamma^{>U}$.
Let $S$ be a subspace of $V$ containing $U$ and let $k_S= \min \{ i \in \{0, \ldots, t \} \mid S \cap V_{i+1} \not= \{0\} \}$ be defined as above.
Then it is easily seen that we have $S\pitchfork_V F$ if and only if $(S+U)\pitchfork_{\overline{V}} \overline{F}$. Moreover, $k_S = \min \{ i \in \{0, \ldots, k_U \} \mid (S+U) \cap (V_{i+1}+U) \not= U \}$. Furthermore, the intersection $S \cap V_{k_S+1}$ is $\omega_{k_S}$-non-degenerate if and only if the intersection $S \cap V_{k_S+1} \cap W$ is $\omega_{k_S}$-non-degenerate. Indeed, in case $k_S < k_U$, there is nothing to show, because $V_{k_S+1} \cap W = V_{k_S+1}$, while in case $k_S = k_U$ the claim follows from the observation that $U \cap V_{k_U+1}$ is $\omega_{k_U}$-non-degenerate and $W$ is $\omega_{k_U}$-orthogonal to $U \cap V_{k_U+1}$. 
The lemma now follows from the fact that, for $x\in S\cap V_{k_S+1} \cap W$, by definition, $x$ is $\omega_{k_S}$-orthogonal to $y \in S\cap V_{k_S+1}\cap W$ if and only if $x+U$ is $\overline{\omega}_{k_S}$-orthogonal to $y+U$. \label{1}
\end{Bw}

The following is an immediate consequence of Lemma \ref{residue}.

\begin{lehilf} \label{residue2}
Let $(\Gamma_j)_j$ be a finite family of $m$ generalized Phan geometries of $V$ and let $U \in \bigcap_j \Gamma_j$. 

Then the geometries $(\bigcap_j \Gamma_j)^{<U} := \{ x \in \bigcap_j \Gamma_j \mid x \subsetneq U \}$ and $(\bigcap_j \Gamma_j)^{>U} := \{ x \in \bigcap_j \Gamma_j \mid x \supsetneq U \}$ are isomorphic to intersections $\bigcap_j \Gamma_j^{<U}$, respectively, $\bigcap_j \Gamma_j^{>U}$ of $m$ generalized Phan geometries of $U$, respectively, $V/U$. In particular, the link of $U$ in $|\bigcap_j \Gamma_j|$ is isomorphic to the join $|(\bigcap_j \Gamma_j^{<U} ) * ( \bigcap \Gamma_j^{>U} )|$ of the intersection of $m$ generalized Phan geometries of $U$ and of the intersection of $m$ generalized Phan geometries of $V/U$.
\end{lehilf}

\begin{Bw}
It is obvious that $(\bigcap_j \Gamma_j)^{<U} = \bigcap_j \Gamma_j^{<U}$ and $(\bigcap_j \Gamma_j)^{>U} = \bigcap_j \Gamma_j^{>U}$. Hence all claims follows from Lemma \ref{residue}.
\end{Bw}

Now we have almost everything in place in order to prove the Main Theorem. It remains to provide the following geometric existence results which help us to determine the homotopy type of a generalized Phan geometry of a two-dimensional vector space. This two-dimensional case will serve as the basis of the induction on the dimension used in the proof of Theorem \ref{allrad} below.

\begin{lehilf} \label{avoid}
Let $m \in \mathbb{N}$, let $\mathbb{F} \leq \mathbb{K}$ be fields, let
$V$ be a vector space of dimension at least two over $\mathbb{F}$, and for $1 \leq i \leq m$ let $\phi_i : V \times V \to \mathbb{K}$ be an $\mathbb{F}$-bilinear map which admits a non-isotropic vector. If $2m < |\mathbb{F}|$, then there exist two $\mathbb{F}$-linearly independent vectors of $V$ which are non-isotropic with respect to all maps $\phi_i$. 
\end{lehilf}

\begin{Bw}
We proceed by induction on the number of maps. For $m=0$ there is nothing to prove. Assuming the theorem is true for $m$ maps, let $v$ be a vector which is non-isotropic with respect to $m$ maps and let $w$ be a vector that is non-isotropic with respect to a $(m+1)$st map. Considering linear combinations of $v$ and $w$ we obtain finitely many quadratic equations in $\mathbb{K}$ whose $\mathbb{F}$-rational zeros have to be avoided in order to find a vector which is non-isotropic with respect to all $m+1$ forms. Since a two-dimensional $\mathbb{F}$-vector space contains $|\mathbb{F}|+1$ one-dimensional subspaces, two vectors of $V$ satisfying the conclusion of the lemma exist, if $2(m+1) < |\mathbb{F}|$. 
\end{Bw}

The same proof works in the setting of hermitian forms of vector spaces over a finite field of square order $q$. The only difference is that there are at most $\sqrt{q}+1$ isotropic one-dimensional subspaces per form in a two-dimensional $\mathbb{F}_q$-vector space instead of two.

\begin{lehilf} \label{avoid2}
Let $V$ be a vector space of dimension at least two over a finite field $\mathbb{F}_{q}$ of square order endowed with $m < \infty$ hermitian forms each of which admits a non-isotropic vector. If $(\sqrt{q}+1)m < q$, then there exist two $\mathbb{F}_{q}$-linearly independent vectors of $V$ which are non-isotropic with respect to all forms.
\end{lehilf}

Finally we can state and prove our main result. Its proof has been strongly influenced by the proof of \cite[Theorem 1.1]{Abels/Abramenko:1993}, \cite[Chapter II, Proposition 12]{Abramenko:1996}.

\begin{thhilf} \label{allrad}
Let $V$ be an $(n+1)$-dimensional vector space over a field $\mathbb{K}$. Let $\sigma$ be an automorphism of $\mathbb{K}$ of order one or two. Let $(\Gamma_j)_{1 \leq j \leq m}$ be a finite family of generalized Phan geometries of $V$. Finally, let $\Gamma = \bigcap_j \Gamma_j$. In case $\mathbb{K} = \mathbb{F}_{q}$ assume $2^nm < q$, if $\sigma$ has order one, and $2^{n-1}(\sqrt{q}+1)m < q$, if $\sigma$ has order two. 

Then $|\Gamma|$ is homotopy equivalent to a (non-trivial) wedge of $(n-1)$-spheres. 
\end{thhilf}

\begin{Bw}
We proceed by induction on $n$. 

For $n=1$, the vector space $V$ is two-dimensional. The geometry $\Gamma = \bigcap_{1 \leq j \leq m} \Gamma_j$ consists of all those one-dimensional subspaces of $V$ that are non-degenerate with respect to each of the $m$ hermitian forms of rank one or two induced by the $\Gamma_j$ on $V$. 
In case $\mathbb{K}$ finite and $\sigma \neq \mathrm{id}$, apply Lemma \ref{avoid2}. Otherwise apply
Lemma \ref{avoid} to the $\mathrm{Fix}_{\mathbb{K}}(\sigma)$-span of a $\mathbb{K}$-basis of $V$. In both cases we find two $\mathbb{K}$-linearly independent vectors of $V$ which are non-isotropic with respect to all forms.
Hence $|\Gamma|$ is a wedge of $0$-spheres.

Assume now that the theorem is true for any family of size $m$ of generalized Phan geometries inside a vector space of dimension $k<n+1$ satisfying $2^{k-1}m < q$, resp.\ $2^{k-2}(\sqrt{q}+1)m < q$ in case $\mathbb{K}$ is a finite field of order $q$.

Let $p$ be a one-dimensional subspace of $V$ which for each of the generalized Phan geometries $\Gamma_j$ is non-degenerate with respect to the unique hermitian form defined on the whole vector space $V$. Such a one-dimensional space $p$ exists by Lemmas \ref{avoid} and \ref{avoid2} via an argument as above.
Define $$Y_0:=\{W\in \Gamma \mid \gen{p,W} \in \Gamma\}.$$ Note that $Y_0$ contains all subspaces of $V$ which contain $p$ and which are contained in $\Gamma$.
Moreover, define $$Y_i=Y_{i-1}\cup\{W\in \Gamma \mid \dim W=n+1-i\}$$ for $i=1,2,\ldots,n$. Obviously, we have $Y_0\subseteq Y_1\subseteq \cdots \subseteq Y_n=\Gamma$. Each set $Y_i$ together with the inclusion relation inside $V$ forms a poset, so that it makes sense to use the symbols $|\cdot|$ and $*$ in this context. 

We are going to use Lemma \ref{mainlemma} and an induction on $i$ to prove that $|Y_i|$ is homotopy equivalent to a wedge of $(n-1)$-spheres for $1\leq i\leq n$.

\begin{description}
\item{$i=0$:} $Y_0$ is contractible, because $U \mapsto \gen{U,p} \mapsto \gen{p}$, $U \in Y_0$, is a chain of deformation retractions.

\item{$i-1 \to i$:}
Apply Lemma \ref{mainlemma} with $Z = |Y_i|$ and $Z' = |Y_{i-1}|$ and $A_j = \overline{St_{Y_i}(U_j)}$ for $U_j \in Y_i \backslash Y_{i-1}$.
The topological space $Z'$ is by induction contractible ($i=1$) or a wedge of $(n-1)$-spheres (otherwise).

To prove $A_{j_1} \cap A_{j_2} \subseteq Z'$, let  $U_{j_1}\neq U_{j_2}\in Y_i\backslash Y_{i-1}$. An element in  $\overline{St_{Y_i}(U_{j_1})}\cap \overline{St_{Y_i}(U_{j_2})}$ is in $Y_i$ and is incident with both $U_{j_1}$ and $U_{j_2}$.  Hence its dimension is not $n+1-i$ and so it lies in $Y_{i-1}$. 

To prove that $A_j$ is contractible
and that $A_j \cap Z'$ is $(n-2)$-spherical, let $U=U_j\in Y_i\backslash Y_{i-1}$. Note that $\dim U= n+1-i$ and $\gen{p,U}\notin \Gamma$, so $p\not\subseteq U$.
Then $$A_j = \overline{St_{Y_{i}}(U)} = |Y_{i-1}^{<U} * U * Y_{i-1}^{>U}|$$ and $$A_j \cap Z' = \overline{St_{Y_{i-1}}(U)} = |Y_{i-1}^{<U} * Y_{i-1}^{>U}|$$ where $Y_{i-1}^{<U}=\{W\in Y_{i-1} \mid W\subsetneq U\}$ and $Y_{i-1}^{>U}=\{W\in Y_{i-1} \mid W \supsetneq U \}$. The topological space $\overline{St_{Y_{i}}(U)} = |Y_{i-1}^{<U} * U * Y_{i-1}^{>U}|$ is homeomorphic to a cone over $|Y_{i-1}^{<U} * Y_{i-1}^{>U}|$, whence contractible.

In order to prove that the topological space $A_j \cap Z' = \overline{St_{Y_{i-1}}(U)} = |Y_{i-1}^{<U} * Y_{i-1}^{>U}|$ is a wedge of spheres, it suffices to prove that both $|Y_{i-1}^{<U}|$ and $|Y_{i-1}^{>U}|$ are wedges of spheres.
As $Y_{i-1}$ contains each element of $\Gamma$ whose dimension is greater than $\mathrm{dim}(U)$, the equality $$Y_{i-1}^{>U}=\Gamma^{>U}=\{W\in \Gamma \mid W \supsetneq U \}$$ holds. By Lemma \ref{residue2} this complex is isomorphic to the intersection of $m$ generalized Phan geometries of $V/U$ (which has dimension $i$), so that by the hypothesis of the induction on $n$ the topological space $|Y_{i-1}^{>U}|$ is a wedge of $(i-2)$-spheres.

As the dimension of each element of $Y_{i-1} \backslash Y_0$ is greater than $\mathrm{dim}(U)$, one easily sees that $$Y_{i-1}^{<U}=Y_0^{<U}=\{W\in \Gamma \mid W\subsetneq U, \gen{W,p}\in \Gamma\}.$$
By Lemma \ref{intersection2} below the latter is isomorphic to the intersection of at most $2m$ generalized Phan geometries of $U$ (which has dimension $n+1-i$). Since, in case $\mathbb{K} = \mathbb{F}_{q}$,  
\begin{eqnarray*}
2^{\dim U-2}C2m & = & 2^{\dim U-1}Cm \\
& \leq & 2^{n-1}Cm \\
& < & q, 
\end{eqnarray*}
where $C$ equals $2$ or $\sqrt{q}+1$ respectively, we can apply the theorem to $Y_{i-1}^{<U}=Y_0^{<U}=\{ W\in \Gamma \mid W\subsetneq U, \gen{W,p}\in \Gamma\}$ by hypothesis of the induction on $n$. Therefore $|Y_{i-1}^{<U}|$ is a wedge of $(n-1-i)$-spheres.

In conclusion  $A_j \cap Z' = \overline{St_{Y_{i-1}}(U)} = |Y_{i-1}^{<U} * Y_{i-1}^{>U}|$ is a wedge of $(n-2)$-spheres, and so by Lemma \ref{mainlemma} the topological space $Z = |Y_i|$ is homotopy equivalent to a wedge of $(n-1)$-spheres.
\end{description}

Consequently, $|\Gamma| = |Y_n|$ is a wedge of $(n-1)$-spheres.
\end{Bw}

A combination of the statement of the theorem with Lemma \ref{residue2} implies the Cohen-Macaulayness of $\Gamma$. 

\begin{kohilf} \label{allrad2}
Let $V$ be an $(n+1)$-dimensional vector space over a field $\mathbb{K}$. Let $\sigma$ be an automorphism of $\mathbb{K}$ of order one or two. Let $(\Gamma_j)_{1 \leq j \leq m}$ be a finite family of generalized Phan geometries of $V$. Finally, let $\Gamma = \bigcap_j \Gamma_j$. In case $\mathbb{K} = \mathbb{F}_{q}$ assume $2^nm < q$, if $\sigma$ has order one, and $2^{n-1}(\sqrt{q}+1)m < q$, if $\sigma$ has order two. 

Then $|\Gamma|$ is Cohen-Macaulay.
\end{kohilf}

\begin{Bw}
Let $s$ be a simplex of $|\Gamma|$. By an iteration of Lemma \ref{residue2} the link $Lk_{|\Gamma|}(s)$ is a join of intersections of generalized Phan geometries. By Theorem \ref{allrad} the link $Lk_{|\Gamma|}(s)$ is homotopy equivalent to a wedge of $\mathrm{dim}(Lk_{|\Gamma|}(s))$-spheres, whence for each $p < \mathrm{dim}(Lk_{|\Gamma|}(s))$ the reduced homology group $\tilde H_p(Lk_{|\Gamma|}(s))$ is trivial. Therefore $|\Gamma|$ is Cohen-Macaulay.
\end{Bw}

We pointed out in Example \ref{examples} that for $t=0$ and $\mathbb{K}$ finite of square order the generalized Phan geometry on a $\mathbb{K}$-vector space equals the complex studied in \cite{Bennett/Shpectorov:2004} to re-prove Phan's first theorem \cite{Phan:1977}. One crucial step in \cite{Bennett/Shpectorov:2004} is the proof of simple connectedness of that geometry. Theorem \ref{allrad} yields more information about the homotopy type of that geometry, however only with an exponential instead of a constant bound on the size of the field. The constant $m$ is replaced by $1$, because we are interested in the homotopy type of one generalized Phan geometry, not the intersection of several generalized Phan geometries. 

\begin{kohilf}
Let $\mathbb{K}$ be a field, let $\sigma$ be an automorphism of $\mathbb{K}$ of order one or two, let $V$ be an $(n+1)$-dimensional vector space over $\mathbb{K}$ endowed with a non-degenerate $\sigma$-hermitian form, and let $\Gamma$ be the complex of proper non-trivial non-degenerate subspaces of $V$. In case $\mathbb{K} = \mathbb{F}_{q}$ assume $2^n < q$, if $\sigma$ has order one, and $2^{n-1}(\sqrt{q}+1) < q$, if $\sigma$ has order two. 

Then $|\Gamma|$ is Cohen-Macaulay. In particular, $|\Gamma|$ is homotopy equivalent to a wedge of $(n-1)$-spheres.
\end{kohilf}

The sphericity of the geometry opposite a chamber is also obtained as a corollary of Theorem \ref{allrad} for a generalized Phan geometry with $t=n$, cf.\ Example \ref{examples2}. As all $\sigma$-hermitian forms $\omega_i$ have rank one and hence the zeros form a hyperplane, we can replace $\sqrt{q}+1$, resp.\ $2$ by $1$ in the bound (cf.\ the proofs of Lemmas \ref{avoid} and \ref{avoid2}) and achieve exactly the same bound as in \cite[Chapter II, Theorem B]{Abramenko:1996}.

\begin{kohilf}
Let $\mathbb{K}$ be a field, let $V$ be an $(n+1)$-dimensional vector space over $\mathbb{K}$, and let $\Gamma$ be the geometry opposite a chamber of $V$. In case $\mathbb{K} = \mathbb{F}_{q}$ assume  $2^{n-1} < q$. 

Then $|\Gamma|$ is Cohen-Macaulay. In particular, $|\Gamma|$ is homotopy equivalent to a wedge of $(n-1)$-spheres.
\end{kohilf}

We still need to provide statement and proof of Lemma \ref{intersection2} which we used for proving Theorem \ref{allrad}.
For this endeavour projections as used in \cite{BGHS2} and \cite{GHS} turn out to be useful. Recall the following fact from linear algebra.

\begin{lehilf} \label{projection}
Let $\mathbb{K}$ be a field, let $\sigma$ be an automorphism of $\mathbb{K}$ of order one or two, let $V$ be an $(n+1)$-dimensional vector space over $\mathbb{K}$, and let $\omega : V \times V \to \mathbb{K}$ be a $\sigma$-hermitian form. Moreover, let $W$ be an arbitrary subspace of $V$ and let
$p$ be a one-dimensional subspace of $V$ which is non-degenerate with respect to $\omega$ and which intersects $W$ trivially.

Then $$\omega^p : V \times V \to \mathbb{K} : (v,w)\mapsto\omega^p(v,w) := \omega\left(\mathrm{pr}_{p^{\perp_{\omega}}}(v),\mathrm{pr}_{p^{\perp_{\omega}}}(w)\right)$$
where $\mathrm{pr}_{p^{\perp_{\omega}}}$ denotes the projection onto the second summand of the direct decomposition $V = p \oplus p^{\perp_{\omega}}$, is a $\sigma$-hermitian form with $$\mathrm{pr}_{W}\left(\Rad\left({\omega}_{|\gen{W,p} \times \gen{W,p}}\right)\right) = \Rad\left({\omega^p}_{|W\times W}\right)$$ where $\mathrm{pr}_{W}$ denotes the projection onto the second summand of the direct decomposition $\gen{W,p} = p \oplus W$.
\end{lehilf}

\begin{Bw}
Let $t\in \Rad({\omega}_{|\gen{W,p}\times\gen{W,p}})$, i.e.\ $t\in \gen{W,p}$ and $\omega(t,x)=0$ for all $x\in \gen{W,p}$, in particular $\omega(t,p)=0$, i.e.\ $t \in p^{\perp_{\omega}}$. Write $t=w+ v$ for $w\in W$ and $v \in p$, so that $\mathrm{pr}_{W}(t)=w$. Notice that, if  $x \in W$, then $\mathrm{pr}_{p^{\perp_{\omega}}}(x)\in \gen{W,p}$. For all $x \in W$ we have $$\omega^p(t,x) = \omega\left(\mathrm{pr}_{p^{\perp_{\omega}}}(t),\mathrm{pr}_{p^{\perp_{\omega}}}(x)\right) = \omega\left(t,\mathrm{pr}_{p^{\perp_{\omega}}}(x)\right)=0,$$ whence $$0=\omega^p(t,x) = \omega^p(w,x) + \omega^p(v,x) = \omega^p(w,x) + \omega(0,\mathrm{pr}_{p^{\perp_{\omega}}}(x)) = \omega^p(w,x)$$ for all $x \in W$, so that $w\in \Rad({\omega^p}_{|W\times W})$.

Conversely, let $w\in \Rad({\omega^p}_{|W\times W})$, i.e.\ $w\in W$ and $$\omega^p(w,x)=\omega\left(\mathrm{pr}_{p^{\perp_{\omega}}}(w),\mathrm{pr}_{p^{\perp_{\omega}}}(x)\right)=0$$ for all $x\in W$. For arbitrary $x \in W$ decomposed as $x = \mathrm{pr}_{p^{\perp_{\omega}}}(x) + v$ for $v \in p$ we have $$\omega\left(\mathrm{pr}_{p^{\perp_{\omega}}}(w),x\right) = \omega\left(\mathrm{pr}_{p^{\perp_{\omega}}}(w),\mathrm{pr}_{p^{\perp_{\omega}}}(x)\right) + \omega\left(\mathrm{pr}_{p^{\perp_{\omega}}}(w),v\right) = \omega\left(\mathrm{pr}_{p^{\perp_{\omega}}}(w),\mathrm{pr}_{p^{\perp_{\omega}}}(x)\right) = 0.$$
On the other hand, $\omega\left(\mathrm{pr}_{p^{\perp_{\omega}}}(w),p\right) = 0$. Hence $\omega\left(\mathrm{pr}_{p^{\perp_{\omega}}}(w),x\right) = 0$ for all $x\in \gen{W,p}$.
 Therefore $\mathrm{pr}_{p^{\perp_{\omega}}}(w) \in \Rad\left(\omega_{|\gen{W,p} \times \gen{W,p}}\right)$, implying $$w = \left(\mathrm{pr}_{W} \circ \mathrm{pr}_{p^{\perp_{\omega}}}\right)(w) \in \mathrm{pr}_{W}\left(\Rad\left(\omega_{|\gen{W,p} \times \gen{W,p}}\right)\right).$$
\end{Bw}

We finish this section with the two final lemmas. Their proofs consist of elementary but extremely technical arguments from linear algebra. 
%We are afraid that no further insight towards the understanding of the proof of Theorem \ref{allrad} can be gained from the remainder of this section except that Theorem \ref{allrad} is correct.

\begin{lehilf} \label{intersection}
Let $\Gamma$ be a generalized Phan geometry of an $(n+1)$-dimensional $\mathbb{K}$-vector space $V$ with respect to a flag $F = \{ V_i \mid 0 \leq i \leq t+1\}$ with forms $(\omega_{i})_{0 \leq i \leq t}$. Let $p$ be an $\omega_t$-non-degenerate one-dimensional subspace of $V$ and let $U \in \Gamma$ be an element of $\Gamma$ not containing $p$. Furthermore, let $$\Delta := \{ x \in \Gamma \mid \gen{x,p} \in \Gamma \}.$$
%
%
%
%
%Moreover, let $\underline{\omega}_i$ be extensions of $\omega_i$ and let $p$ be a one-dimensional subspace of $V \backslash V_t$ as in the statement of Lemma \ref{extension}, and let $U$ be an element of $\Gamma$ not containing $p$. Furthermore, let $\Delta$ be the subgeometry of $\Gamma$ consisting of all those elements of $\Gamma$ whose span with $p$ is also contained in $\Gamma$. 
If $\Gamma^{<U} = \{ x \in \Gamma \mid x \subsetneq U \}$ is non-empty, then $\Gamma^{<U} \cap \Delta$ is an intersection of at most two generalized Phan geometries of $U$. 
\end{lehilf}

\begin{Bw}
The geometry $\Gamma^{<U}$ is isomorphic to the generalized Phan geometry of $U$ with respect to the flag $\{V_i \cap U \mid k_U \leq i \leq t+1\}$ and forms $\left({\omega_{i}}_{|(U \cap V_{i+1}) \times (U \cap V_{i+1})}\right)_{k_U \leq i \leq t}$, cf.\ Lemma \ref{residue} and its proof. 
To determine the isomorphism type of the intersection $\Gamma^{<U} \cap \Delta$, let $W \lneq U$ be an element of $\Gamma^{<U}$. The vector space $W$ is contained in $\Gamma^{<U} \cap \Delta$ if and only if $\gen{p,W}$ is an element of $\Gamma$. By definition, we have $\gen{p,W} \in \Gamma$ if and only if $\gen{p,W} \pitchfork_V F$  with $\gen{p,W} \cap V_{k_{\gen{p,W}+1}}$ non-degenerate with respect to $\omega_{k_{\gen{p,W}}}$ where, as before, $k_{\gen{p,W}}: = \min \{ i \in \{0, \ldots, t \} \mid \gen{p,W} \cap V_{i+1} \not= \{0\} \}$.

Clearly, $\gen{p,W}\pitchfork_V F$ if and only if $W\pitchfork_U F^p_U$ where $F^p_U := \{\gen{V_i,p} \cap U \mid 0 \leq i \leq l_U \}$ with $l_U := \min\{0 \leq i \leq t+1 \mid \gen{V_{i},p} \cap U = U \}$. 
Note that $l_U$ is always equal to $t+1$, except when $V_t$ is a hyperplane of $V$, in which case $l_U=t$.

Using notations
 from Lemma \ref{projection}, we define 
$$\underline{\omega}_i^p : V \times V \to \mathbb{K} : (v,w)\mapsto \underline{\omega}_i\left(\mathrm{pr}_{p^{\perp_{\underline{\omega}_{i}}}}(v),\mathrm{pr}_{p^{\perp_{\underline{\omega}_{i}}}}(w)\right)$$ for $0 \leq i \leq l_U$.

We have $k_{\gen{p,W}}=\min \{ i \in \{0, \ldots, l_U \} \mid W \cap \gen{V_i,p} \cap U \not= \{0\} \}=:j$.
Moreover, we claim that $\gen{p,W} \cap V_{j+1}$ is non-degenerate with respect to $\omega_{j}$  if and only if $W \cap (\gen{V_{j+1},p}\cap U)$ is non-degenerate with respect to ${{\underline\omega}_{j}^p}_{|(\gen{V_{j+1},p} \cap U) \times (\gen{V_{j+1},p} \cap U)}$.
The intersection $\gen{p,W} \cap V_{j+1}$ is non-degenerate with respect to $\omega_{j}$ if and only if $\gen{\gen{p,W} \cap V_{j+1}, p}= \gen{W,p} \cap \gen{V_{j+1},p}=\gen{W \cap \gen{V_{j+1},p},p}$ is non-degenerate with respect to $\underline{\omega}_j$ (this is obvious for $j=t$ and follows from the fact that $p \perp_{\underline{\omega}_j} V_{j+1}$ by \ Lemma \ref{extension} for $j<t$).
Since $p$ is non-degenerate with respect to ${\underline{\omega}_j}$ and $p\not\subseteq W \cap \gen{V_{j+1},p}$, we can apply Lemma \ref{projection} and obtain $$\mathrm{pr}_{W \cap \gen{V_{j+1},p}}\left(\Rad\left({\underline{\omega}_j}_{|\gen{W \cap \gen{V_{j+1},p},p} \times \gen{W \cap \gen{V_{j+1},p},p}}\right)\right) = \Rad\left({\underline{\omega}_j^p}_{|(W \cap \gen{V_{j+1},p}) \times (W \cap \gen{V_{j+1},p})}\right).$$
Therefore $\gen{W \cap \gen{V_{j+1},p},p}$ is non-degenerate with respect to $\underline{\omega}_j$ if and only if $\gen{V_{j+1},p}\cap W$  is non-degenerate with respect to $\underline{\omega}_j^p$, because $p$ is non-degenerate with respect to ${\underline{\omega}_j}$. Since  $\gen{V_{j+1},p}\cap W\subset U$, the above claim is proved.

Notice that the family $(\gen{V_i,p} \cap U)_{0 \leq i \leq l_U}$ may contain some subspace of $U$ more than once. Since $\dim(\gen{V_i,p} \cap U) = \dim(V_i \cap U) + 1$ if and only if $p \subseteq \gen{V_i,U}$ and $i\neq t+1$, and since the subspaces $V_i$ form a flag, the only subspaces of $U$ that can occur more than once in that family are zero-dimensional or one-dimensional. 
Notice that there is at least one one-dimensional subspace in $F_U^p$ exactly when $p \subseteq \gen{V_{k_U},U}$.
Since the $\sigma$-hermitian form on the one-dimensional subspace of the flag of a generalized Phan geometry does not really play any role | none of the transversal spaces have a non-trivial intersection with this space | we can choose an arbitrary non-degenerate $\sigma$-hermitian form $\omega$ on that subspace (which will be the first form of the generalized Phan geometry).    

It remains to check the radical property $$\Rad\left({{\underline{\omega}_i^p}}_{|(\gen{V_{i+1},p} \cap U) \times (\gen{V_{i+1},p} \cap U)}\right) = \gen{V_i,p} \cap U$$ for the relevant incides $i$. Let us start with $k_U + 1 \leq i \leq t$.
By Lemma \ref{projection}, we have $$\mathrm{pr}_{\gen{V_{i+1},p} \cap U}\left(\Rad\left({\underline{\omega}_i}_{|\gen{\gen{V_{i+1},p} \cap U,p} \times \gen{\gen{V_{i+1},p} \cap U,p}}\right)\right) = \Rad\left({\underline{\omega}_i^p}_{|(\gen{V_{i+1},p} \cap U) \times (\gen{V_{i+1},p} \cap U)}\right).$$
Since $\gen{\gen{V_{i+1},p} \cap U,p} = \gen{V_{i+1},p} \cap \gen{U,p}$, it suffices to determine the $\underline{\omega}_i$-radical of $\gen{V_{i+1},p} \cap \gen{U,p}$ | which is 
$V_i \cap \gen{U,p}$ | and to project it onto $\gen{V_{i+1},p} \cap U$ | which results in $\gen{V_i,p} \cap U$. 

In view of our comments on the spaces $V_i$ with $\gen{V_i,p} \cap U$ one-dimensional, the only remaining relevant index is $i = k_U$. In case $\gen{V_{k_U},p} \cap U$ one-dimensional, we deduce 
that $p$ is contained in the span of $V_{k_U}$ and $U$, so that by symmetry $V_{k_U} \cap \gen{U,p}$ is one-dimensional as well. Since $V_{k_U+1}\cap \gen{U,p}=(U\cap V_{k_U+1})\oplus (V_{k_U} \cap \gen{U,p})$, this implies that $V_{k_U} \cap \gen{U,p}$ is the $\underline{\omega}_{k_U}$-radical of $\gen{V_{k_U+1},p} \cap \gen{U,p}$ and, thus, its projection $\gen{V_{k_U},p} \cap U$ onto $\gen{V_{k_U+1},p} \cap U$ is the $\underline{\omega}_{k_U}^p$-radical, so that the arguments from above apply. 
In this case, $F_U^p$ is in fact $\{\{0\}\}\cup\{\gen{V_i,p} \cap U \mid k_U\leq i \leq l_U\}$, to which we associate the forms $\left(\omega, \left({\underline{\omega}_{i}^p}_{|(\gen{V_{i+1},p} \cap U) \times (\gen{V_{i+1},p} \cap U)}\right)_{k_U\leq i \leq l_U} \right)$. We just showed that they satisfy the radical property.

In case $\gen{V_{k_U},p} \cap U$ zero-dimensional, however, the space $\gen{V_{k_U+1},p} \cap \gen{U,p}$ may be $\underline{\omega}_{k_U}$-non-degenerate or may have a one-dimensional $\underline{\omega}_{k_U}$-radical, and so $\gen{V_{k_U+1},p} \cap U$ may be ${\underline{\omega}_{k_U}^p}$-non-degenerate or may have a one-dimensional ${\underline{\omega}_{k_U}^p}$-radical $R$.  

In the first case $R=\Rad\left({\underline{\omega}_{k_U}^p}_{|(\gen{V_{k_U+1},p} \cap U) \times (\gen{V_{k_U+1},p} \cap U)}\right)=\{0\}=\gen{V_{k_U},p} \cap U$, and  $F_U^p$ is in fact $\{\{0\}\}\cup\{\gen{V_i,p} \cap U \mid k_U+1 \leq i \leq l_U\}$, to which we associate the forms $\left(\omega, \left({\underline{\omega}_{i}^p}_{|(\gen{V_{i+1},p} \cap U) \times (\gen{V_{i+1},p} \cap U)}\right)_{k_U+1\leq i \leq l_U} \right)$. We just showed that they satisfy the radical property.

In the second case, notice that $F^p_U$ does not contain a one-dimensional space and so we can define a slightly modified flag ${F^p_U}^+ := F^p_U \cup \left\{ R\right\}$, which is in fact $\{\{0\},R\}\cup\{\gen{V_i,p} \cap U \mid k_U+1 \leq i \leq l_U\}$. Define also an arbitrary non-degenerate $\sigma$-hermitian form $$\omega_+ : R \times R \to \mathbb{K}.$$ We then associate to ${F^p_U}^+$ the forms $\left(\omega_+, \left({\underline{\omega}_{i}^p}_{|(\gen{V_{i+1},p} \cap U) \times (\gen{V_{i+1},p} \cap U)}\right)_{k_U+1\leq i \leq l_U} \right)$. Now  ${F^p_U}^+$ and these forms satisfy the radical property.
Since we changed the flag, there are still a few things to be checked in order to show that $\gen{p,W} \in \Gamma$ if and only if $W$ is in the generalized Phan geometry of $U$ with respect to the flag ${F^p_U}^+$ and forms as described. 
Assume $\gen{p,W} \in \Gamma$. We already proved that $W$ is transversal to $F_U^p$, so we just need to show that $W\cap R=\{0\}$. This intersection is contained in $\Rad\left({\underline{\omega}_{k_U}^p}_{|(\gen{V_{k_U+1},p} \cap W) \times (\gen{V_{k_U+1},p} \cap W)}\right)$. If $\gen{V_{k_U+1},p} \cap W=\{0\}$ we are done. If $\gen{V_{k_U+1},p} \cap W\neq\{0\}$, then $\gen{W,p} \cap V_{k_U+1}\neq\{0\}$, and so $k_{\gen{p,W}}=k_U$ (since $\gen{W,p} \cap V_{k_U}=\{0\}$). Therefore by hypothesis $\gen{W,p} \cap V_{k_U+1}$ is non-degenerate with respect to $\omega_{k_U}$, and so, as shown above, $W \cap \gen{V_{k_U+1},p}$ is non-degenerate with respect to ${{\underline\omega}_{k_U}^p}_{|(\gen{V_{k_U+1},p} \cap U) \times (\gen{V_{k_U+1},p} \cap U)}$. Hence $W\cap R=\{0\}$. Finally, we have that $k_{\gen{p,W}}=\min \{ i \in \{k_U+1, \ldots, l_U \} \mid W \cap (\gen{V_i,p} \cap U) \not= \{0\} \}$, and so we can conclude by the same argument as above. Since $W\pitchfork_U {F^p_U}^+$ implies $W\pitchfork_U {F^p_U}$ and $W\cap R$ must be trivial, the reverse implication is immediate.

In conclusion, the intersection $\Gamma^{<U} \cap \Delta$ is isomorphic to the intersection of the generalized Phan geometry of $U$ with respect to the flag $(V_i \cap U)_{k_U \leq i \leq t+1}$ and forms $\left({\omega_{i}}_{|(V_{i+1} \cap U) \times (V_{i+1} \cap U)}\right)_{k_U \leq i \leq t}$ and the generalized Phan geometry of $U$ with respect to the flag $F^p_U$ (respectively, ${F^p_U}^+$) and forms as described above.
\end{Bw}

Similar to Lemma \ref{residue2} the following is immediate from Lemma \ref{intersection}.

\begin{lehilf} \label{intersection2}
Let $(\Gamma_j)_j$ be a finite family of $m$ generalized Phan geometries of $V$ 
with respect to flags $F^j = \{ V^j_i \mid 0 \leq i \leq t_j+1\}$ with forms $(\omega^j_{i})_{0 \leq i \leq t_j}$. Let $p$ be a one-dimensional subspace of $V$ which for each $j$ is $\omega^j_{t_j}$-non-degenerate. Let $U \in \bigcap_j \Gamma_j$ be an element of $\bigcap_j \Gamma_j$ not containing $p$. Furthermore, let $$\Delta := \{ x \in \bigcap_j \Gamma_j \mid \gen{x,p} \in \bigcap_j \Gamma_j \}.$$
If each $\Gamma_j^{<U} = \{ x \in \Gamma_j \mid x \subsetneq U \}$ is non-empty, then $\bigcap_j \Gamma_j^{<U} \cap \Delta$ is an intersection of at most $2m$ generalized Phan geometries of $U$. 
\end{lehilf}

\begin{Bw}
This follows from Lemma \ref{intersection}, because $(\bigcap_j \Gamma_j^{<U}) \cap \Delta = \bigcap_j (\Gamma_j^{<U} \cap \Delta_j)$ where $\Delta_j:= \{ x \in \Gamma_j \mid \gen{x,p} \in \Gamma_j \}$.
\end{Bw}

\section{Outline of applications} \label{outline}

In this section we sketch two sets of applications of Theorem \ref{allrad} 
to geometric group theory. A thorough treatment of 
these applications requires the language of (twin) buildings which we do
not want to develop here. For more details and background information
we refer the reader to \cite{DM} and \cite{Gramlich}. 

Our applications depend on the following fact. 
Example \ref{examples3} indicates how to prove its correctness 
for a (twin) building $\mathcal{B}$ of type $A_n$. A special case is treated in \cite{Blok/Hoffman}.

\begin{fact} \label{fact}
Let $\mathcal{B}$ be a twin building with a Phan involution $\tau$. Let $R$ be a spherical residue of type $\bigoplus_{i \in I} A_{n_i}$ of $\mathcal{B}$. Then the geometry associated to the chamber system $A_\tau(R)$ is a direct sum of generalized Phan geometries of type $A_{n_i}$.
\end{fact}

A Phan involution $\tau$ is an involution interchanging both parts of a twin building isometrically and satisfying an additional technical condition. The chamber system $A_\tau(R)$ consists of all those chambers of $R$ with minimal possible codistance from their image under $\tau$, cf.\ Example \ref{examples3}. Detailed discussions can be found in \cite{DM} and \cite{Gramlich}.

\medskip
The first set of applications of Theorem \ref{allrad} deals with the revision program of Phan's results \cite{Phan:1977}, \cite{Phan:1977a}.

\begin{thhilf} \label{thm}
Let $\mathbb{K}$ be a field containing sufficiently many elements, let $\sigma$ be an automorphism of $\mathbb{K}$ of order one or two, let $\Delta$ be a graph without triangles, and let $\mathcal{B}$ be a twin building of type $\Delta$ over $\mathbb{K}$ with a Phan involution $\tau$ that induces $\sigma$ on $\mathbb{K}$.

Then the Phan geometry $\mathcal{B}_\tau$, i.e., the subsystem of chambers of the positive half of $\mathcal{B}$ mapped onto an opposite one by $\tau$, is $2$-simply connected. 
\end{thhilf}

\begin{sketch}
By the criterion of simple connectedness from \cite{DM} it remains to prove that for residues $R$ of rank two, resp.\ three the chamber system $A_\tau(R)$ is connected, resp.\ simply connected. Since therefore one only has to deal with the connectedness properties of chambers systems of rank two or three, by \cite{Tits:1981} one can equivalently study connectedness properties of the corresponding geometries. By Fact \ref{fact} it suffices to study generalized Phan geometries of type $A_2$ and $A_3$ or direct sums of generalized Phan geometries. By Theorem \ref{allrad} the generalized Phan geometries of type $A_2$ are connected and the generalized Phan geometries of type $A_3$ are simply connected, if the field $\mathbb{K}$ contains sufficiently many elements. A direct sum of two rank one geometries is always connected and a direct sum of a rank one and a connected rank two geometry is always simply connected. Hence the claim follows.
\end{sketch}

The following result is implied by Theorem \ref{thm}. 
It has first been stated in \cite{Phan:1977a} 
and then re-proven via a direct computation 
based on \cite{DM} by Hoffman, Shpectorov and two of the authors during
an RiP stay in Oberwolfach in summer 2005. For the definition of a group with a weak Phan system, 
we refer the reader to the 
survey \cite{Bennett/Gramlich/Hoffman/Shpectorov:2003}, 
the thesis \cite{Gramlich:2004}, 
and the research paper \cite{Bennett/Shpectorov:2004}.

\begin{kohilf} \label{Phan}
Let $n \geq 3$, let $q \geq 5$, let $\Delta \in \{ A_n, D_n, E_6, E_7, E_8 \}$, and let $K$ be a group with
a weak Phan system of type $\Delta$ over $\mathbb{F}_{q^2}$. Then $K$ is a central quotient of $\mathrm{SU}_{n+1}(q^2)$, $\mathrm{Spin}^\pm_{2n}(q)$, or the simply connected version of ${}^2E_6(q^2)$, $E_7(q)$, $E_8(q)$, respectively.
\end{kohilf}

\begin{sketch}
Since $q \geq 5$, the inequalities $2(q+1)<q^2$ and $2^2(q+1) < q^2$ are satisfied, so that by Theorem \ref{allrad} and Theorem \ref{thm} the Phan geometry $\mathcal{B}_\tau$ with respect to the building of type $\Delta$ over $\mathbb{F}_{q^2}$ and the unique semilinear Phan involution $\tau$, is simply connected. By \cite{Bennett/Gramlich/Hoffman/Shpectorov:2003}, \cite{Bennett/Shpectorov:2004}, \cite{Gramlich:2004} the simple connectedness of $\mathcal{B}_\tau$ implies that $K$ is the universal enveloping group of the Phan amalgam contained in $K$. The claim follows from the classification of unambiguous noncollapsing Phan amalgams in \cite{Bennett/Shpectorov:2004}, \cite{Gramlich:2004}, by which there exists a unique Phan amalgam up to passing to quotients.
\end{sketch}

As for the second set of applications the Main Theorem can be used in order to establish topological finiteness properties of so-called unitary forms of Kac-Moody groups over finite fields.  
The following statement is a first rough attempt.

\begin{kohilf} \label{finitely generated}
Let $n \geq 3$, let $q \geq 5$, let $\Delta$ be a graph without triangles, and let $G$ be a Kac-Moody group of type $\Delta$ over $\mathbb{F}_{q^2}$. Let $\tau$ be the product of the Chevalley involution of $G$ and the field involution. 

Then the group $K := \mathrm{Fix}_G(\tau)$ is finitely presented.
\end{kohilf}

\begin{sketch}
As in the sketch of Corollary \ref{Phan}, the Phan geometry $\mathcal{B}_\tau$ with respect to the twin building of type $\Delta$ over $\mathbb{F}_{q^2}$ is simply connected. Again as before the simple connectedness of $\mathcal{B}_\tau$ implies that $K$ is the universal enveloping group of the Phan amalgam contained in $K$. Since this Phan amalgam consists of finite groups, the group $K$ is finitely presented.
\end{sketch}

More sophisticated statements are possible if one uses the existing machinery for topological finiteness properties, notably Brown's criteria \cite{Brown:1987}. The following theorem deals with the topological finiteness length of the unitary form of a Kac-Moody group of type $\tilde A_n$ over the field $\mathbb{F}_{q^2}$.

\begin{thhilf} \label{inv}
Let $n \geq 2$, let $q$ be be a prime power satisfying $2^{n-1}(q+1) < q^2$, and let $\sigma$ be the involution 
of the field $\mathbb{F}_{q^2}(t)$ mapping $t$ to $t^{-1}$ and acting 
as the field involution on $\mathbb{F}_{q^2}$. Then the group 
$$K := \mathrm{SU}_{n+1}(\mathbb{F}_{q^2}[t,t^{-1}], \sigma) = 
\{ A \in \mathrm{SL}_{n+1}(\mathbb{F}_{q^2}[t,t^{-1}]) 
\mid (A^{-1})^\sigma = A^T \}$$ is  of type $F_{n-1}$ and not of type $F_n$.
\end{thhilf}

\begin{sketch}
This follows by Brown's criterion \cite[Corollary 3.3]{Brown:1987} applied to the 
filtration $(C_n)_{n \in \mathbb{N}}$ from \cite{DM} of the positive half $\mathcal{B}_+$ of the twin building of type $\widetilde{A}_n$ 
for the Kac-Moody group $\mathrm{SL}_{n+1}(\mathbb{F}_{q^2}[t,t^{-1}])$. 
The relative links again are generalized Phan geometries by Fact \ref{fact}, 
which are spherical by Theorem \ref{allrad}. The group $K$ is a subgroup of the locally compact group $\mathrm{SL}_{n+1}(\mathbb{F}_{q^2}((t)))$ and since $\sigma$ interchanges $t$ and $t^{-1}$ it is necessarily discrete. As the cell stabilisers of the action of $\mathrm{SL}_{n+1}(\mathbb{F}_{q^2}((t)))$ on $\mathcal{B}_+$ are compact, the cell stabilisers in its discrete subgroup $K$ have to be finite. Therefore the cell stabilisers in $K$ are of type $F_\infty$. Finally, using the surjectivity of the norm map $\mathbb{F}_{q^2} \to \mathbb{F}_q$ one can prove that for each $n \in \mathbb{N}$ the group $K$ acts transitively on the layer $C_n \backslash C_{n-1}$ of the filtration, so that the filtration is cocompact. 
\end{sketch}

\begin{bem}
It is elementary but not straightforward to see that the group $\mathrm{SU}_{n+1}(\mathbb{F}_{q^2}[t,t^{-1}], \sigma)$ from Theorem \ref{inv} is an arithmetic subgroup of $\mathrm{SL}_{n+1}(\mathbb{F}_{q^2}(t))$. Hence, for the sake of brevity of this remark, we allow ourselves the liberty to use some very heavy machinery in order to establish this arithmeticity:
By \cite{GM} the group $\mathrm{SU}_{n+1}(\mathbb{F}_{q^2}[t,t^{-1}], \sigma)$ is a lattice in the locally compact group $\mathrm{SL}_{n+1}(\mathbb{F}_{q^2}((t)))$. Therefore by \cite[Chapter IX]{Margulis:1991} it is an arithmetic subgroup of $\mathrm{SL}_{n+1}(\mathbb{F}_{q^2}(t))$. This means that Theorem \ref{inv} states that the topological finiteness length of $\mathrm{SU}_{n+1}(\mathbb{F}_{q^2}[t,t^{-1}], \sigma)$ equals its local rank minus one, which is another piece of evidence in favour of the rank conjecture on the finiteness lengths of $S$-arithmetic groups; cf.\ \cite{Behr:1998}.
\end{bem}

%The involution described in Theorem \ref{inv} is an example of a Phan involution of the group $\mathrm{SL}_{n+1}(\mathbb{F}_{q^2}[t,t^{-1}])$. We refer the reader to \cite{Blok/Hoffman} for a discussion of alternative Phan involutions of this group.

\end{document}